\documentclass[12pt,thmsa]{article}
\usepackage{amsfonts}
\usepackage{amssymb}
\newtheorem{teo}{Theorem}[section]

\newtheorem{obs2}[teo]{Remark}

\newtheorem{tea}{Theorem}[subsection]

\newtheorem{no2}[teo]{Note}

\newtheorem{no3}[tea]{Note}

\newcommand{\Q}{{\mathbb{Q} }}

\newcommand{\GL}{{\rm GL}}

\newcommand{\GSp}{{\rm GSp}}

\title{Galois characterization of Endoscopy for rational Siegel modular forms}
\author{Luis V. Dieulefait\\ Centre de Recerca Matem{\'a}tica\\ Bellaterra,
Barcelona
  \thanks{
e-mail: LDieulefait@crm.es.
Supported by MECD postdoctoral fellowship at the Centre de Recerca Matem{\'a}tica from
the Ministerio de Educaci{\'o}n y Cultura
}} 

\begin{document}

\maketitle
\begin{abstract} We establish a relation between Galois reducibility and Endoscopy for
 genus 2 Siegel cusp forms which have rational eigenvalues and are  unramified at $3$.

\end{abstract}

\section{The Theorem}
Let $f$ be a genus $2$ Siegel cuspidal Hecke eigenform of weight
$k >2$ and $\Q_f$ the number field generated by its eigenvalues.
It is well known that if $f$ is not of Saito-Kurokawa type but it
is  ``endoscopic" (i.e., it is in the image of the weak endoscopic
lift) the compatible family of Galois representations $ \{ \rho_{
f, \lambda} \} $ attached to $f$ (constructed by Taylor, Laumon
and Weissauer for any Siegel cusp form) will be reducible
over $\Q_f$ , with two $2$-dimensional irreducible components. \\
In this note we will prove that the converse statement is true,
for the case $\Q_f = \Q$. We will have to impose a local condition
at $3$ and will assume that the determinants are minimally
ramified. More precisely, the result is:
\begin{teo}
  \label{teo:endos}
  Let $f$ be a genus $2$ Siegel cusp form of weight $k>2$ with
  $\Q_f = \Q$, such that the corresponding automorphic
  representation $\pi_f$ has multiplicity one and $\pi_{f,3 }$, its local
  component at $3$, is unramified. Assume that the compatible
  family of Galois representations $\{ \rho_{ f,
\ell} \}$ attached to $f$ reduces (over $\Q$) as follows:
$$ \rho_{ f,
\ell} \simeq \sigma_{1,\ell}  \oplus \sigma_{2,\ell}  \qquad
\qquad (1.1)$$ for every prime $\ell$, where $\{  \sigma_{1,\ell}
\}$  and $\{  \sigma_{2,\ell} \}$ are compatible families of
$2$-dimensional irreducible representations both with
determinant $\chi^{2k-3}$.\\
Then, $f$ is endoscopic. More precisely, there exist two classical
cuspidal modular forms $f_1$, $f_2$, of weights $2$ and $2k-2$
(respectively) such that the family of representations $ \{
\sigma_{1,\ell} \otimes \chi^{2-k} \}$ is attached to $f_1$ and
the family $\{ \sigma_{2,\ell} \} $ is attached to $f_2$.

  \end{teo}
Remarks: 1- The irreducibility assumption of the $2$-dimensional
components is equivalent to assume that $f$ is not of
Saito-Kurokawa type.\\
2- Reducibility of the whole family $\{ \rho_{ f, \ell} \}$  as in
the statement of the theorem is equivalent to a similar condition
imposed only at a single prime $\ell$, provided that $\ell > 4k-5$
and $\pi_{f,\ell}$ is unramified. This follows from a result of
``existence of a family" proved in [D2].\\

\section{Proof}

Irreducibility of the $2$-dimensional components implies that we
are not in the Saito-Kurokawa case, and together with the
multiplicity one assumption this implies (as proved by Weissauer)
that the representations are pure, odd, and for every prime $\ell$
such that $\pi_{f,\ell}$ is unramified, the representations
$\sigma_{1,\ell}$ and $\sigma_{2,\ell}$ are crystalline with
Hodge-Tate weights $\{ k-2, k-1 \}$ and $\{ 0, 2k-3 \}$,
respectively. \\
Let us first show modularity of the family $\{ \sigma_{1,\ell}
\otimes \chi^{2-k} \}$. By assumption, $\pi_{f,3}$ is unramified,
thus $\sigma_{1,3} \otimes \chi^{2-k} $ is a Barsotti-Tate
representation, irreducible, odd, with rational coefficients, and
unramified outside a finite set of primes. Then, applying a
combination of modularity results of Diamond-Taylor-Wiles and
Skinner-Wiles (as done in [D1] and [D2]) we conclude that
$\sigma_{1,3} \otimes \chi^{2-k} $ is modular, and this gives
modularity of the family $ \{ \sigma_{1,\ell} \otimes \chi^{2-k}
\}$. The corresponding modular form $f_1$ must have weight $2$
because for almost every $\ell$ the representations in this family
are
Barsotti-Tate.\\
This argument ``{\`a} la Wiles" can not be applied to $\sigma_{2,3}$
because, even if we  again have Wiles' starting point (namely, we
know that residually $ \bar{\sigma}_{2,3}$ is either modular or
reducible), the prime $3$ is too small compared with the
difference $2k-3$ of the Hodge-Tate weights to make the
strategy workable.\\
To show modularity of the family $\{  \sigma_{2,\ell} \}$  we will
explode the fact that (1.1) is telling us that the representations
$\sigma_{2,\ell}$ can be obtained by ``substracting" a modular
representation from another modular representation.\\
A key ingredient is a result recently proved by Weselmann (yet
unpublished, but see [BWW] and [W]), which states that $\pi_f$ has
a weak lift to an automorphic representation $\pi'$ of $\GL(4,
\mathbb{A})$, where $\mathbb{A} $ are the rational adeles. Thus,
by Cebotarev, the family $\{  \rho_{f,\ell} \}$ is also attached
to $\pi'$. \\
We want to apply a result of Jacquet and Shalika (which appears as
theorem 3.3 in [T2]), in a similar way than what is done in [T2],
section 533. We have from (1.1) the equality of $L$-functions:
$$ L(\pi' , s) = L(\sigma_{2,\ell} , s) L( f_1 \otimes \chi^{k-2} , s) $$
Observe that $\pi'$ is the weak lift of $\pi_f$, but it is
not necessarily cuspidal.\\
To conclude that $\sigma_{2,\ell}$ is modular, as in section 533
of [T2], we must find a prime $\ell$ such that $L(
\sigma^*_{2,\ell} \otimes \sigma_{2,\ell} , s )$ has a simple pole
at $s=1$, because in that case the result of  Jacquet and Shalika
implies $\sigma_{2,\ell} \simeq \sigma_{ \pi_i ,\ell}$, where
$\pi_i$ is one of the cuspidal constituents of $\pi'$. Then, it
only remains to find a prime satisfying this condition.\\
Take $\ell > 4k-5$ such that the local components of $\pi'$ and
$\pi_f$ at $\ell$ are unramified. For such a prime $\ell$ the
representation $\sigma_{2,\ell}$ is crystalline with Hodge-Tate
weights $\{0, 2k-3 \}$ and the main result of [T1] implies that
there exists a totally real number field $F$ such that the
restriction of $\sigma_{2,\ell}$ to the Galois group of $F$ is
modular, i.e., it agrees with the Galois representation attached
to a Hilbert modular form over $F$.\\
But, as explained in [T2], section 533, precisely from this
potentially modular property (and the fact that it is preserved
after solvable base change) one can deduce that $L(
\sigma^*_{2,\ell} \otimes \sigma_{2,\ell} , s )$ does have a
simple pole at $s=1$, as we wanted. This shows modularity of the
family $\{ \sigma_{2,\ell} \}$ and it is clear from its Hodge-Tate
decomposition that it corresponds to a modular form of weight
$2k-2$. We conclude that the Siegel cusp form $f$ is endoscopic.

\section{Bibliography}


[BWW] Ballmann, J., Weissauer, R., Weselmann, U., {\it Remarks on
the Fundamental Lemma for stable twisted Endoscopy of Classical
Groups}, preprint (2002)\\
\newline
[D1] Dieulefait, L., {\it Modularity of Abelian Surfaces with Quaternionic
Multiplication}, Math. Res. Letters {\bf 10} no. 2-3 (2003)\\
\newline
[D2] Dieulefait, L., {\it Existence of families of Galois
representations and new cases of the Fontaine-Mazur conjecture},
preprint (2003)\\
available at http://www.arxiv.org/math.NT\\
\newline
%
[T1] Taylor, R., {\it On the meromorphic continuation of degree two
 L-functions}, preprint (2001)\\
 available at http://abel.math.harvard.edu/$\sim$rtaylor/ \\
\newline
[T2] Taylor, R., {\it Galois Representations}, proceedings of ICM,
Beijing, August 2002\\
longer version available at
http://abel.math.harvard.edu/$\sim$rtaylor/ \\
\newline
[W] Weissauer, R., {\it A remark on the existence of Whittaker
models for $L$-packets of automorphic representations of
$\GSp(4)$}, preprint


\end{document}